\def \NN {\mathbb N}

\def \RR {\mathbb R}

\def \epsilon{\varepsilon}

\def \S  {{\mathcal S}}

\def \si {\sigma}

\newcommand{\res}{{\rm res}}
\renewcommand{\S}{{\mathcal S}}

\documentclass[12pt,reqno]{amsart}

\usepackage{amsmath,amssymb}

\newtheorem{lem}{Lemma}

\usepackage{url} 
\usepackage{xcolor}

\begin{document}
\title[]{On the invariants of $L$-functions of degree $2$, II: \\ twisted conductor}
\author[]{J.Kaczorowski \lowercase{and} A.Perelli}
\maketitle

\hfill {\it Dedicated to Helmut Maier on the occasion of his 70th birthday}

\medskip
{\bf Abstract.} {\small We continue our study of the invariants of the $L$-functions of deree 2 twisted by Dirichlet characters. In a previous paper, we proved that twisting does not change degree and internal shift. Our aim now is to show a corresponding theorem for the conductor. This paper is intended as the second part of a series aiming to prove a general version of the Weil converse theorem with minimal assumptions on the shape of the functional equation of twists.}



\section{Introduction}

We continue our study of the invariants of the $L$-functions of deree 2 twisted by Dirichlet characters. We keep the general notation introduced in \cite{Ka-Pe/InvariantsOfTwistedI}, where we proved that twisting does not change degree and internal shift. Here we show a corresponding result for the conductor. This paper is intended as the second part of a series aiming to prove a general version of the Weil converse theorem with minimal assumptions on the shape of the functional equation of twists. For the readers' convenience, we recall all the necessary notation in the next section. 

\medskip
{\bf Theorem.} {\sl Let $F\in {\mathcal S}^{\sharp}$ be of degree $2$ with integer conductor $q_F$. Let $D$ be a square-free positive integer, $D\equiv 1\, (\bmod \, q_F)$, such that $F$ splits polynomially of degree $2$ at every prime $p|D$. Moreover, suppose that $F$ is saturated $(\bmod \, D)$ and for every Dirichlet character $\chi \, (\bmod \, D)$, the twist $F^{\chi^*}$ belongs to $\S^{\sharp}$. Then the conductors of $F^{\chi^*}$, $F$ and  $\chi$ are related by}
\begin{equation}
\label{twistchar}
q_{F^{\chi^*}}= q_F f_{\chi}^2.
\end{equation}

\medskip
We expect that a suitable version of \eqref{twistchar}, namely
\begin{equation}
\label{twistchar2}
q_{F^\chi} = q_F f_\chi^{d_F},
\end{equation}
holds in much greater generality. In particular, we expect \eqref{twistchar2} to hold for every $F\in\S$, with any degree $d_F$ and integer conductor $q_F$, and every primitive character $\chi$\,(mod\,$f_\chi$), $(f_\chi,q_F)=1$, such that the twist $F^\chi$ also belongs to $\S$; see conjecture (III.4) in \cite{Ka-Pe/2000}. Note that \eqref{twistchar2} holds for all classical $L$-functions.

 \medskip
 {\bf Acknowledgements}.  This research was partially supported by the Istituto Nazionale di Alta Matematica, by grant PRIN 2022 {\sl ``The arithmetic of motives and L-functions''} and by grant 2021/41/BST1/00241 {\sl ``Analytic methods in number theory''}  from the National Science Centre, Poland.

\section{Notation}

In this paper 
\[
F(s)=\sum_{n=1}^{\infty}\frac{a(n)}{n^s}\ \ \ \ \ (\sigma>1)
\]
stands for a generic function from the extended Selberg class $\S^{\sharp}$. Recall that every such function satisfies a functional equation of the form
\[
\gamma(s) F(s)=\omega \overline{\gamma(1-\overline{s})} \overline{F(1-\overline{s})},
\]
where
\[
\gamma(s)= Q^s\prod_{j=1}^r\Gamma(\lambda_j s+\mu_j)
\]
and $Q>0$, $\lambda_j >0$, $\Re(\mu_j)\geq 0$, $ |\omega|=1$. The Selberg class $\S$ is the subclass of the functions $F\in\S^\sharp$ having a rather general Euler product expansion and satisfying the Ramanujan conjecture $a(n)\ll n^\varepsilon$. Degree, conductor, internal shift and root number of $F$ are defined, respectively, by
\[
d_F=2\sum_{j=1}^r \lambda_j,\quad q_F=(2\pi Q)^2\prod_{j=1}^r\lambda_j^{2\lambda_j},\quad \theta_F= \Im(\sum_{j=1}^r\mu_j),\quad \omega^*_F= \omega\prod_{j=1}^r\lambda_j^{-2i\Im(\mu_j)}.
\]
We refer to our survey papers \cite{Ka-Pe/1999b},\cite{Kac/2006},\cite{Per/2005},\cite{Per/2004} for definitions, examples, and the basic theory of the Selberg class.

\smallskip
For a Dirichlet character $\chi$ and real numbers $\alpha$ and $\beta$, we denote by $F^{\chi}(s)$ and $F(s,\beta,\alpha)$ the twists
\[
F^{\chi}(s)= \sum_{n=1}^{\infty} \frac{a(n)\chi(n)}{n^s} \quad \text{and} \quad F(s,\beta,\alpha)= \sum_{n=1}^{\infty} \frac{a(n)}{n^s}e\left(-\beta n-\alpha\sqrt{n}\right).
\]
If $F^\chi\in\S^\sharp$, then its degree, conductor, internal shift and root number are denoted, respectively, by $d_{F^\chi}$, $q_{F^\chi}$, $\theta_{F^\chi}$ and $\omega_{F^\chi}^*$. In the important special case $\beta=0$ and $\alpha\neq0$, $F(s,0,\alpha)$ is called the standard twist of $F$ and is denoted by $F(s,\alpha)$. Concerning the standard twist, we define the spectrum of $F$ as
\[
\text{Spec}(F) = \{\alpha>0: a(n_\alpha)\neq 0\}, \ \text{where} \ n_\alpha = q_Fd_F^{-d_F}\alpha^{d_F} \  \text{and} \  a(x) = 0 \  \text{if} \ x\not\in\NN.
\]
The relevance of Spec$(F)$ to the polar structure of $F(s,\alpha)$ is reported in Lemma \ref{lem:3}.

\smallskip
We say that $F\in \S^{\sharp}$ splits polynomially of degree 2 at a prime $p$ if 
\[ 
F(s) = \prod_{j=1}^2\left(1-\frac{\alpha_j(p)}{p^s}\right)^{-1} \sum_{p\nmid n} \frac{a(n)}{n^s}
\]
for certain complex numbers $\alpha_j(p)\neq 0$. We say that $F\in \S^{\sharp}$ is saturated $(\bmod \, D)$ if for every integer $M$ and every residue class $a\, (\bmod \,D)$, $(a,D)=1$, there exists a positive integer $\nu\equiv a\,(\bmod \,D)$ with $(\nu,M)=1$ such that $a(\nu)\neq 0$. Moreover, we say that $F\in \S^{\sharp}$ is normalized if $\theta_F=0$ and the first non-zero Dirichlet coefficient of $F$ equals $1$.

\smallskip
For a Dirichlet character $\chi$, we denote by $\chi^*$ the primitive character inducing $\chi$. As usual, we denote by $f_{\chi}$ the conductor of $\chi$, thus $\chi^*$ is a Dirichlet character $(\bmod \,f_{\chi})$. Sometimes, in cases where no confusion may occur, when $\chi$ is a character $(\bmod \,q)$ we simply write $q^*$ instead of $f_{\chi}$. Finally, the Gauss sum of $\chi$ is denoted by $\tau(\chi)$.

\section{Lemmas}

\begin{lem}
\label{lem:1} 
Let $D$ be a positive squarefree integer and $(a,D)=1$, and suppose that $F\in \mathcal{S}^{\sharp}$ with $d_F=2$ splits polynomially of degree $2$ at all primes $p|D$. Then for $\si>1$ we have
\begin{equation}
\label{eq:fchi1}
F(s,a/D,\alpha) = \frac{1}{\varphi(D)} \sum_{\chi(\bmod D)} \sum_{m|\left(\frac{D}{f_{\chi}}\right)^2}
\frac{B_m(F)f(\chi,m,a)}{m^s} F^{\chi^*}(s,\sqrt{m}\alpha),
\end{equation}
where
\[
f(\chi,m,a)= \mu\left(\frac{D}{f_{\chi}}\right) \overline{\tau(\chi^*)}\overline{\chi^*\left(\frac{D}{f_{\chi}}\right)}
\chi^*(am) \prod_{p|m}p \neq 0.
\]
\end{lem}
\begin{proof}
This is a special case of Lemma 3 in \cite{Ka-Pe/InvariantsOfTwistedI} with $\lambda=1/2$. Note that in such a Lemma 3, the analog of equation \eqref{eq:fchi1} has $m|(\frac{D}{f_{\chi}})^*$ in place of $m|(\frac{D}{f_{\chi}})^2$. Here we have simply $m|(\frac{D}{f_{\chi}})^2$ thanks to the hypothesis that $F$ splits polynomially of degree $2$ at all primes $p|D$.
\hfill
\end{proof}
\begin{lem}
\label{lem:2}  
Let $F\in\S^\sharp$ with $d_F=2$ be normalized. Then for  $\beta>0$ and $\alpha\in\RR$ we have
\[
F(s,\beta,\alpha) =-\omega_F^* e^{-i\frac{\pi}{2}\xi_F}  e\big(\frac{\alpha^2}{4\beta}\big) (\beta\sqrt{q_F})^{2s-1}  \sum_{n=1}^\infty \frac{\overline{a(n)}}{n^s} e\left(\frac{n}{\beta q_F} - \frac{\alpha}{\beta\sqrt{q_F}}\sqrt{n}\right) + h(s)
\]
for $\sigma>1$, where $h(s)$ is holomorphic for $\si>1/2$.
\end{lem}
\begin{proof}
This is the Lemma in \cite{Ka-Pe/2017} with a slightly changed notation.
\hfill
\end{proof}
\begin{lem}
\label{lem:3}
The standard twist $F(s,\alpha)$ of  $F\in \S^{\sharp}$ is meromorphic on $\mathbb C$ for every $\alpha\in\RR\setminus\{0\}$. Moreover, if $F\in \S^{\sharp}$ has $d_F=2$ and $\theta_F=0$, then for $\alpha\in$ {\rm Spec}$(F)$ we have
\begin{equation}
\label{residue}
\res_{s=3/4} F(s,\alpha)= c_0(F) \frac{ \overline{a(n_{\alpha})}}{n_{\alpha}^{1/4}},
\end{equation}
where $c_0(F)\neq 0$ is a constant depending on $F$ and $n_{\alpha}=q_F\alpha^2/4$.
\end{lem}
\begin{proof}
This follows from Theorems 1 and 2 in \cite{Ka-Pe/2005}.
\hfill
\end{proof}
We remark that \eqref{residue} holds also if $0<\alpha\not\in$ Spec$(F)$, in the sense that both sides equal zero in this case.
\begin{lem}
\label{lem:4}
Let $F$, $D$ and $\chi$ be as in the Theorem. Then $d_{F^{\chi^*}}=2$ and $\theta_{F^{\chi^*}}=\theta_F$.
\end{lem}
\begin{proof}
This is the main result of part I of the series, see \cite{Ka-Pe/InvariantsOfTwistedI}. 
\hfill
\end{proof}
\begin{lem}
\label{lem:5}
If the Theorem is true for every entire $F\in \S^{\sharp}$ with $d_F=2$, then it holds true for all $F\in\S^{\sharp}$ with $d_F=2$.
\end{lem}
\begin{proof}
Let $F\in \S^{\sharp}$ with $d_F=2$ have a pole of order $m_F>0$ at $s=1$, fix a prime $p_0\nmid D$ and write
\[ 
G(s)=\left(1-\frac{p_0+1}{p_0^s}+ \frac{p_0}{p_0^{2s}}\right)^{m_F}.
\]
We easily check that $G\in\S^{\sharp}$, and has degree $0$ and conductor $q_G=p_0^{2m_F}$. Thus the function $F_1(s)=G(s)F(s)$ belongs to $\S^{\sharp}$, is entire and has $d_{F_1}=2$ and $q_{F_1}= q_Gq_F= p_0^{2m_F}q_F$.
Moreover, it is easily seen that $G^\chi\in\S^\sharp$ for every character $\chi$, and $q_{G^\chi}=q_G$. Hence applying the Theorem to $F_1$ we get
 \[ 
 q_{F_1^{\chi^*}}= q_{F_1} f_{\chi}^2.
 \]
 But $F_1^{\chi}= G^{\chi}F^{\chi}$ and $q_{F_1^{\chi^*}}= q_{G^{\chi^*}}q_{F^{\chi^*}}= p_0^{2m_F}q_{F^{\chi^*}}$. Therefore
 \[ 
 p_0^{2m_F}q_{F^{\chi^*}} = p_0^{2m_F}q_F f_{\chi}^2
 \]
 and the result follows.
 \hfill
 \end{proof}
\begin{lem}
\label{lem:3a}
Let $\psi(\bmod \,q)$ be a Dirichlet character induced by $\psi^*(\bmod \, q^*)$. Then
\[
\tau(\psi)=\mu\left(q/q^*\right) \psi^*(q/q^*)\tau(\psi^*).
\]
\end{lem}
\begin{proof}
This is a well-known property of the Gauss sums, see Theorem 9.10 in \cite{M-V}.
\end{proof}

\section{Proof of the Theorem}

Using Lemma \ref{lem:5}, we may assume without loss of generality that $F$ is entire. Shifting the variable $s$ by a suitable purely imaginary constant and multiplying by an appropriate non-zero constant, we can further assume that  $F$ is normalized, thus we may apply Lemma \ref{lem:2}. Taking into account that $D\equiv 1(\bmod \,q_F)$, for $\si>1$ we obtain
\begin{eqnarray*}
F(s,1/D,\alpha)&=& e^{a_1+b_1s}e\left(\frac{\alpha^2D}{4}\right)
\sum_{n=1}^{\infty}\frac{\overline{a(n)}}{n^s}e\left(\frac{nD}{q_F} -\frac{\alpha D}{\sqrt{q_F}}\sqrt{n}\right) +h(s)\\
&=&e^{a_1+b_1s}e\left(\frac{\alpha^2D}{4}\right)
\sum_{n=1}^{\infty}\frac{\overline{a(n)}}{n^s}e\left(\frac{n}{q_F} -\frac{\alpha D}{\sqrt{q_F}}\sqrt{n}\right) +h(s)\\
&=&e^{a_1+b_1s}e\left(\frac{\alpha^2D}{4}\right)
\overline{\sum_{n=1}^{\infty}\frac{a(n)}{n^{\overline{s}}}e\left(\frac{nD}{q_F} -\frac{\alpha D}{\sqrt{q_F}}\sqrt{n}\right)} +h(s),
\end{eqnarray*}
where $h(s)$ is holomorphic for $\si>1/2$. Applying Lemma \ref{lem:2} once more we see that the last expression equals
\begin{equation}
\label{eq:F1D-2}
\begin{split}
e^{a_2+b_2s}e\left(\frac{\alpha^2D}{4}(1-D)\right)
&\overline{\sum_{n=1}^{\infty}\frac{\overline{a(n)}}{n^{\overline{s}}}e\left(\alpha D\sqrt{n}\right)} +h(s) \\ 
&= e^{a_2+b_2s}e\left(\frac{\alpha^2D}{4}(1-D)\right) F(s,\alpha D) + h(s).
\end{split}
\end{equation}
Next, applying Lemma \ref{lem:1} with $a=1$, we see that for $\sigma>1$
\begin{equation}
\label{eq:F1D-1}
F(s,1/D,\alpha) = \frac{1}{\varphi(D)} \sum_{\chi(\bmod D)} \sum_{m|\left(\frac{D}{f_{\chi}}\right)^2}
\frac{B_m(F)f(\chi,m,1)}{m^s} F^{\chi^*}(s,\sqrt{m}\alpha).
\end{equation}
Comparing (\ref{eq:F1D-1}) and (\ref{eq:F1D-2}), we arrive at the following identity involving various standard twists 
\begin{equation}
\label{eq:F=SumF}
\begin{split}
e^{a_2+b_2s}e\left(\frac{\alpha^2D}{4}(1-D)\right)&F(s,\alpha D) + h(s)  \\
&=\frac{1}{\varphi(D)} \sum_{\chi(\bmod D)} \sum_{m|\left(\frac{D}{f_{\chi}}\right)^2}\frac{B_m(F)f(\chi,m,1)}{m^s} F^{\chi^*}(s,\sqrt{m}\alpha).
\end{split}
\end{equation}
Note that \eqref{eq:F=SumF} is proved for $\sigma>1$, but both sides are meromorphic for $\sigma>1/2$ thanks to Lemmas \ref{lem:2} and \ref{lem:3}; thus \eqref{eq:F=SumF} holds in this larger half-plane. Observe that $\theta_{F^{\chi^*}}=0$ by Lemma \ref{lem:4} and hence, according to Lemma \ref{lem:3}, $F^{\chi^*}(s,\sqrt{m}\alpha)$ has a pole at $s=3/4$ whenever $\sqrt{m}\alpha\in$ Spec$(F^{\chi^*})$. Taking residues of both sides of (\ref{eq:F=SumF}) and using Lemma \ref{lem:3} again, we obtain
\begin{equation}
\label{eq:res}
\begin{split}
A_0(F)& e\left(\frac{\alpha^2D}{4}(1-D)\right)
\frac{\overline{a\left(q_F\alpha^2D^2/4\right)}}{\left(q_F\alpha^2D^2/4\right)^{1/4}} \\
&=\frac{1}{\varphi(D)} \sum_{\chi(\bmod D)} \sum_{m|\left(\frac{D}{f_{\chi}}\right)^2}\frac{B_m(F)f(\chi,m,1)}{m^{3/4}}
c_0(F^{\chi^*})\frac{\overline{a(q_{\chi}m\alpha^2/4)}\overline{\chi(q_{\chi}m\alpha^2/4)}}{(q_{\chi}m\alpha^2/4)^{1/4}},
\end{split}
\end{equation}
where $A_0(F)\neq 0$ and, for simplicity, we write $q_{\chi}$ instead of $q_{F^{\chi^*}}$.

\smallskip
Let 
\[ 
k(\chi,m):= \frac{q_{\chi}}{q_F} \frac{m}{D^2}\ \ \ \ \text{where}\ \ \ \ m\mid  \left(\frac{D}{f_{\chi}}\right)^2,
\]
\[
k(\chi):= k(\chi, (D/f_{\chi})^2),\quad k_0:=\max_{\chi(\bmod D)} k(\chi).
\]
Now in (\ref{eq:res}) we choose
\[
\alpha = \frac{2}{D}\sqrt{\frac{\nu}{k_0q_F}}\qquad (\nu\in {\mathbb N}),
\]
hence
\[ 
\frac{\alpha^2}{4}D(1-D) = \frac{\nu}{k_0D}\frac{1-D}{q_F},\qquad  \frac{q_F}{4}\alpha^2D^2 = \frac{\nu}{k_0}, \qquad \frac{q_{\chi}m}{4} \alpha^2 = \frac{q_{\chi}}{q_F}\frac{m}{D} = \frac{k(\chi,m)}{k_0}.
\]
Thus (\ref{eq:res}) implies that
 \begin{equation}
 \label{eq:res2}
 e\left(\frac{\nu}{k_0D}\frac{1-D}{q_F}\right)\overline{a\left(\frac{\nu}{k_0}\right)} =
  \sum_{\chi(\bmod D)} \sum_{m|\left(\frac{D}{f_{\chi}}\right)^2}c_F(\chi,m)
  \overline{a\left(\frac{k(\chi,m)\nu}{k_0}\right)} \overline{\chi(\nu)}  
 \end{equation}
with certain complex  coefficients $c_F(\chi,m)$ such that 
\begin{equation}
\label{eq:c(chi)}
c_F(\chi, (D/f_{\chi})^2)\neq 0.
\end{equation} 

\smallskip
Next we show that $k_0=1$. We observe that $k_0 \geq k(\chi_0,D^2) = 1$ since $f_{\chi_0}=1$.  Suppose, by contradiction, that $k_0>1$. If $k_0\in {\mathbb Q}$ then we write 
\[ 
k_0 = \frac{a}{b}\ \ \ \ (a,b)=1 , a>1.
\]
Let 
\[ 
a^*:= \begin{cases} 
a & \text{ if \ $k_0\in {\mathbb Q}$}\\
1 & \text{ if \ $k_0\not\in {\mathbb Q}$}.
\end{cases}
\]
Therefore, we have that $\nu/k_0\not\in {\mathbb N}$ if $(\nu,a^*)=1$, and hence $a\left({\nu}/{k_0}\right)=0.$ Recalling (\ref{eq:res2}), for all $(\nu,a^*)=1$ we have that
\begin{equation}
\label{eq:res3}
\mathop{\sum_{\chi(\bmod D)}\sum_{m\mid (D/f_{\chi})^2}}_{k(\chi,m)/k_0\in {\mathbb Q}}
\overline{c_F(\chi,m)} a\left(\frac{k(\chi,m)\nu}{k_0}\right) \chi(\nu)=0.
\end{equation}
When the quotient $k(\chi,m)/k_0$ is rational we write
\[ 
\frac{k(\chi,m)}{k_0} = \frac{A(\chi,m)}{B(\chi,m)}
\]
with $A(\chi,m),B(\chi,m)\in {\mathbb N}$ and $(A(\chi,m),B(\chi,m))=1.$ If $k(\chi,m)<k_0$ then $B(\chi,m)>1$. Let now
\[ 
M:=  \mathop{\prod_{\chi(\bmod D)}\prod_{m\mid(D/f_{\chi})^2}}_{k(\chi,m)/k_0\in {\mathbb Q}, k(\chi,m)<k_0} B(\chi,m).
\]
For $(\nu,M)=1$ we have 
\[ 
\frac{k(\chi,m)}{k_0}\nu = \frac{A(\chi,m) \nu}{B(\chi,m)}\not\in {\mathbb N},
\]
and hence for such $\nu$ we see that
\[ 
a\left(\frac{k(\chi,m)}{k_0}\nu\right) =0.
\]
Moreover, for $k(\chi,m)=k_0$ we have $m=(D/f_{\chi})^2$. Thus, for $(\nu,M)=1$,  (\ref{eq:res3}) reads as
\[ 
a(\nu)\sum_{\substack{\chi(\bmod D)\\ k(\chi)=k_0}} c_F(\chi) \chi(\nu)=0,
\]
where $c_F(\chi):= c(\chi, (D/f_{\chi})^2)$. Note that, according to (\ref{eq:c(chi)}), these coefficients are non-zero.
Since $F$ is saturated $(\bmod \,D)$, for every resudue class $a\, (\bmod \,D)$, $(a,D)=1$, we can find $\nu\equiv a\,(\bmod \,D)$, $(\nu, M)=1$, such that $a(\nu)\neq 0$. Hence for every $(a,D)=1$ we have
\[
\sum_{\substack{\chi(\bmod D)\\ k(\chi)=k_0}} c_F(\chi) \chi(a)=0.
\]
This contradicts the linear independence of characters, thus we have $k_0=1$. 

\smallskip
Finally we show that $k(\chi)=1$ for every character $\chi\,(\bmod \,D)$. Since $k_0=1$, (\ref{eq:res2}) reads as
\begin{equation}
\label{formula} 
 e\left(\frac{\nu}{D}\frac{1-D}{q_F}\right)\overline{a(\nu)} =
 \mathop{ \sum_{\chi(\bmod D)} \sum_{m|\left({D}/{f_{\chi}}\right)^2}}_{k(\chi,m)\in {\mathbb Q}}
c_F(\chi,m) \overline{a(k(\chi,m)\nu)} \overline{\chi(\nu)}. 
\end{equation}
For $(\nu, M)=1$, we have $k(\chi, m)\nu\not\in {\mathbb N}$, and hence the corresponding terms on the right hand side vanish. In particular, non-zero terms must have $k(\chi,m)=1$. Recalling that in the latter case $m=(D/f_{\chi})^2$, equation \eqref{formula} takes the form
\[ 
 e\left(\frac{\nu}{D}\frac{1-D}{q_F}\right)\overline{a(\nu)} =  \overline{a(\nu)}
  \sum_{\substack{\chi(\bmod D)\\k(\chi)=1}}c_F(\chi)\overline{\chi(\nu)}. 
\]
As before, we use the assumption that $F$ is saturated $(\bmod \,D)$, thus getting that
\begin{equation}
\label{eq:schi}
 e\left(\frac{a}{D}\frac{1-D}{q_F}\right)= \sum_{\substack{\chi(\bmod D)\\k(\chi)=1}}c_F(\chi)\overline{\chi(a)}
  \end{equation}
for all $(a,D)=1$. Suppose that the sum on the right hand side of \eqref{eq:schi} does not contain all characters $(\bmod \,D)$. For instance, suppose that $k(\psi)\neq 1$ for a certain $\psi(\bmod \,D)$. Multiplying both sides of (\ref{eq:schi}) by $\psi(a)$ and summing over all $(a,D)=1$ we obtain, by the orthogonality of characters, that
\[
\sum_{\substack{a(\bmod D)\\(a,D)=1}}\psi(a) e\left(\frac{a}{D}\frac{1-D}{q_F}\right)=0.
\]
Recalling Lemma \ref{lem:3a}, we see that this sum equals
\[
\overline{\psi\left(\frac{1-D}{q_F}\right)}\tau(\psi) =\overline{\psi\left(\frac{1-D}{q_F}\right)}\mu(D/f_{\chi}) \psi^*(D/f_{\chi}) \tau(\psi^*)
\] 
and hence is non-zero since $D$ is square-free. This contradiction shows that $k(\chi)=1$ and consequently $q_{\chi}=q_Ff_{\chi}^2$ for all $\chi\,(\bmod \,D)$. The proof is complete.

\bigskip
\noindent
Jerzy Kaczorowski, Faculty of Mathematics and Computer Science, A.Mickiewicz University, 61-614 Pozna\'n, Poland. e-mail: \url{kjerzy@amu.edu.pl}

\medskip
\noindent
Alberto Perelli, Dipartimento di Matematica, Universit\`a di Genova, via Dodecaneso 35, 16146 Genova, Italy. e-mail: \url{alberto.perelli@unige.it}

\end{document}